\documentclass[12pt]{amsproc}
\everymath{\displaystyle}
\usepackage{amsmath,amsfonts,amsthm}
 at 10truept
\font\regsmallcaps=cmcsc12 at 12truept
\def\c{\mathbb{C}}
\def\higherpowers{{\rm\ higher\ powers}}
\def\be{\begin{enumerate}}
\def\ee{\end{enumerate}}
\def\bi{\begin{itemize}}
\def\ei{\end{itemize}}
\def\ms{\medskip}
\def\bs{\bigskip}
\def\bc{\begin{center}}
\def\ec{\end{center}}
\def\schroeder{{Schr\"{o}der}}
\def\equaldef{ \buildrel\hbox{\scriptsize  def}\over= }
\begin{document}
\title{Elements of Finite Order in the Group of
Formal  Power Series Under Composition}
\author{Marshall M. Cohen}
\address{Department of Mathematics, Morgan State University,
1400 E. Cold Spring Lane, Baltimore, MD 21251}
\email{marshall.cohen@morgan.edu}
\newtheorem{theorem}{Theorem}
\newtheorem*{unnumberedtheorem}{Theorem}
\newtheorem*{corollary*}{Corollary}
\newtheorem{lemma}{Lemma}[section]
\newtheorem{corollary}[lemma]{Corollary}
\newtheorem{notation}[lemma]{Notation}
\newtheorem{example}[lemma]{Example}
\newtheorem{remark}[lemma]{Remark}
\newtheorem{definition}[lemma]{Definition}
\begin{abstract}
We consider formal power series \ $f(z) = \omega z + a_2z^2 + \ldots \ (\omega \neq 0)$, \quad
with coefficients in a field of characteristic $0$. These form a group under the operation of composition
(= substitution). We prove (Theorem 1) that every element $f(z)$ of finite order is conjugate to its linear term
$\ell_\omega(z) = \omega z$, and we characterize those elements which conjugate $f(z)$ to $\omega z$. Then we investigate the construction of elements of order $n$ and prove (Theorem 2) that, given
a primitive $n$'th root of unity $\omega$ and an arbitrary sequence $\{a_k\}_{k\neq nj+1}$  there is a unique
sequence $\{a_{nj + 1}\}_{j=1}^\infty$ such that the series \ $f(z) = \omega z +
a_2z^2 + a_3z^3 + \ldots$\ \ \  has order $n$. Sections 1 - 5 give an exposition of this classical subject, written for the 2005 - 2006 Morgan State University Combinatorics Seminar. We do not claim priority for these results in this classical field, though perhaps the proof of Theorem 2 is new. We have now (2018) added Section 6 which gives references to valuable articles in the literature and historical comments which, however incomplete,  we hope will give proper credit to those who have preceded this note and be helpful and of interest to the reader.
\end{abstract}
\maketitle

\section{Introduction}
Suppose that  $F$ is a field  of characteristic 0. We denote
by \  $F[[z]]$\ the set of all formal power series with coefficients
in $F$.  Composition in $F[[z]]$ is defined for elements of the form\quad
$f(z) = a_0 + a_1z + a_2z^2 +
\cdots,\quad  g(z) =  b_1z + b_2z^2 + \cdots $\quad  by:
$$ (f\circ g)(z) \equiv f(g((z)) =a_0 + a_1g(z) + a_2g(z)^2 + \ldots
\ = \ \sum_{n \geq 0}\ a_ng(z)^n. $$
{\bf Define:}\quad $G = G[[z]] \equiv \{\,f \in F[[z]]\ \Big\vert \ f(z) = \sum_n a_nz^n, \quad
a_0 = 0,\  a_1 \neq 0\,\}$.

It is a classical fact (See Section 6.1) that $G$ is a group under
composition. 
The identity element of $G$ is the series  $id$ given by $id(z) = z$. We
denote
the inverse of $f$ by $\overline{f}$ and the composition \ $f\circ f\circ f \cdots \circ f\ (n\  {\rm times})$ \ by \ $f^{(n)}$. We define $\ell_\omega \in G$ by $\ell_\omega(z) =
\omega z$.

\medskip
The purpose of this note is to prove the following two theorems concerning formal power series of finite compositional order. In Section 5, we comment on series of order two and series of infinite order.

\begin{theorem}[Classification of finite order elements up to
 conjugation]\mbox{}\\
Suppose that \ $f(z) = \omega z + \sum_{k=2}^{\infty}a_kz^k$ is an element of  $G$ of order $n$.\ \  Define $f^\ast$ by
$$f^\ast = \frac{1}{n}\,\sum_{j=1}^n\, \omega^{n-j}f^{(j)}.$$
 Then:
 
 \ms
 \be
 \item [A.] $f^\ast \in G,\ \ n$ is the multiplicative order of\  $\omega \in F - \{0\}$\ and
 \[ f^\ast\circ f\circ\overline{f^\ast}\,(z) = \omega z.\]

\ms
{\noindent Two elements of order $n$ in $G$ are conjugate iff their
lead coefficients are the same primitive $n\text{'}$th root of unity.}

\bs
\item [B.] 
\begin{enumerate}
\item If $g \in G$ then
\[\qquad g\circ f\circ \overline{g} = \ell_\omega\  \iff \ \exists
 h(z) = \sum_{j=0}^\infty h_{nj+1}z^{nj+1}\ \in G\
 \mbox{such that} \ g = h\circ f^\ast.\]
\item For every sequence \ $\{\,g_{nj+1}\}_{j=0}^{\infty}$\
 of
elements of $F$ with $g_1\neq 0$ there exists a unique sequence \ $\{g_k\}_{1 < k \neq nj + 1}$ \  such that
the series $g(z) = \sum_{k=1}^\infty g_kz^k$ satisfies
\ $g\circ f\circ \overline{g} = \ell_\omega$.
\ee%
\ee
\end{theorem}
\bigskip
\noindent {\bf Note:} \quad  The same proof of Theorem 1A, in which the conjugating element $f^\ast$ is used, is  given by Cheon and Kim [5] and by O'Farrell and Short [10]. (See Section 6.) 
\ms
When $F = \c$, the field of complex numbers, Theorem 1A. is a very special case of Theorem 8 (p. 19) of [1]. The proof given here, using $f^{\ast}$, is our interpretation of page 23 of [1] in the case of a finite cyclic group of formal power series of a single variable. The definition of $f^{\ast}$ is a finite version of the formal power series S in equation (85) on page 23 of [1]. See also the final 
Remark in Section 5 of this paper.

\medskip
If $F = \c$  and $f$ is analytic (has a positive radius of convergence) of compositional order $n$ then
clearly $f^\ast$ is also analytic. Thus $f$ is conjugate by  an
analytic function in a neighborhood of the origin  to a rotation
of order $n$. (See Theorem 5 (p. 55) of [1], for a more general result.)

\medskip
In the light of Theorem 1 one there are uncountably many different series of compositional order $n$  which have a given primitive $n$'th root of unity $\omega$ as lead  coefficient --- namely, the conjugates of $\ell_\omega$.
How much freedom is there in choosing the other coefficients?

\bigskip
\begin{theorem}
[{The coefficients of finite order elements of
$G[[z]]$}]\mbox{}

If\ $n$\ is a positive integer and $\omega$ is a primitive $n'{\rm
th}$\ root of unity in the field $F$ then for every infinite sequence  \
$\{a_k\}_{1 < k \neq nj + 1}$ \ of elements of $F$ there exists a unique sequence \ $\{\,a_{nj+1}\}_{j=1}^{\infty}$
\ such that the formal power series

 \centerline{$f(z) = \omega z + \sum_{k=2}^{\infty}\ a_k\,z^k$}

 has order $n$ in $G[[z]]$.
\end{theorem}
\begin{corollary*}
If $F = \mathbb{C}$ and $2 \leq n \in \mathbb{N}$ then there exist uncountably many elements $f(z) \in G[[z]]$ of order $n$ such that the power series $f(z)$ does not converge for any non-zero $z\in \mathbb{C}$.
\end{corollary*}
\bigskip

{\bf Acknowledgments:}\quad I am grateful to A. Nkwanta for
introducing me to these matters and to the other members of the
2005 - 2006 Combinatorics Seminar at Morgan State University, B. Eke, X. Gan
and L. Woodson for their interest and encouragement. I would
like to thank H. Furstenberg and L. Shapiro for their helpful comments.
\bigskip

\section{The n{\it th} power of a series under composition}
\medskip
We do a number of elementary calculations based on the fact that
if $g(z) = b_1z + b_2z^2 + \cdots$\ then the multiplicative (as
opposed to compositional) power\  $g(z)^k$\ has as its lead term
\ $b_1^kz^k$:
$$g(z)^k =  (b_1z + b_2z^2 + \cdots)(b_1z + b_2z^2 + \cdots) \ldots
(b_1z + b_2z^2 + \cdots) = b_1^kz^k + \higherpowers.$$
To start with, this combines with the definition of composition to give
\medskip
\begin{lemma}

If\ $f(z) =  a_1z + a_2z^2 +
\cdots$,\ \ and \ \   $g(z) =  b_1z + b_2z^2 + \cdots$ then
$$(f\circ g)(z) = (a_1b_1)z + (a_1b_2 + a_2b_1^2)z^2
+ (a_1b_3 +2a_2b_1b_2 + a_3b_1^3)z^3 +
\higherpowers.\quad \Box$$
\end{lemma}
\medskip
{\bf Notation:} We write\quad $f^{(n)}(z) = \sum_{k=1}^\infty\,f^{(n)}_k\,z^k$.
\bigskip
The following well known lemma can be proved by induction.
We omit the proof since it is a corollary to Lemma 2.4 below.
The fact that $F$ is a field of characteristic 0 is used.
\medskip
\begin{lemma}  If $f(z) =  z + a_kz^k + a_{k + 1}z^{k + 1} + \cdots$,\ with  $a_k \neq 0$ and $n$  a positive integer then
\[f^{(n)}(z) = (z + a_knz^k + \higherpowers)\]
and $f$ has infinite order in $G$.\quad $\Box$.
\end{lemma}
\bigskip
\begin{lemma} If \ \ $f(z) = \omega z + a_2z^2 + \cdots\ \in G$\ \
has finite order  under composition then the order of $f$ equals the
order of $\omega \in F - \{0\}$.
\end{lemma}
\medskip
{\bf Proof:}\quad Let order$(f) = n$. Then \ $f^{(n)}(z) = z = \omega^nz$ by repeated use of Lemma 2.1. Thus
$\omega^n = 1$. If $a = {\rm order}(\omega)$\ let \ $n = ab$.   Then
\begin{eqnarray*}
f^{(a)}(z) &=& \omega^az + \higherpowers = z + \higherpowers.\\
            f^{(n)}(z)  &=& z = (f^{(a)})^{(b)}(z).
\end{eqnarray*}

By Lemma 2.2, we see that $f^{(a)}(z) = z$.
Therefore \ $a = {\rm order}(\omega) = n = {\rm order}(f).\ \Box$

\bigskip
 Suppose that\  $f(z) =a_1z + a_2z^2 + \ldots\ \in G$.
 From the definition of $f^{(n)}_k$ and Lemma 2.1 we see that \[f^{(1)}_k = a_k,\qquad  f^{(n)}_1 = a_1^n,\qquad f^{(2)}_2 = a_2a_1(1 + a_1))\qquad f^{(2)}_3 = a_3a_1(1 + a_1^2) +2a_1a_2^2\]
More generally
we have the following key lemma, which will be used in proving the uniqueness in Theorem 2.
\begin{lemma} If $n, k$ are positive integers then there exists a polynomial \ $P^{(n)}_k(x_1, \ldots, x_{k-1})$ \  in k-1 variables, with integer coefficients, such that
\ $P^{(n)}_k(x_1, 0, 0, \ldots,\, 0) = 0$ (i.e., each summand
contains an $x_j$ with $j\neq 1$),\ and such that $P^{(n)}_k$ satisfies the following:
\medskip
\mbox{} \qquad If $f(z) = (\omega z + a_2z^2 + \cdots a_{k }z^{k} + \cdots\ ) \in G[[z]]$ then
$$f^{(n)}_k  = a_k\omega^{n - 1}\bigg(1 + \omega^{k - 1} +
\cdots\  \omega^{(k-1)(n-1)}\bigg) +  P^{(n)}_k(\omega, a_2, \ldots, a_{k-1})$$
\end{lemma}
\medskip
{\bf Proof:} We proceed by induction on $n$

\medskip
$n = 1:\quad f^{(1)}_k = a_k = a_k\omega^0(1)$ \ and\ $P^{(1)}_k=0$.
\medskip
Suppose $n > 1$. For simplicity we introduce the further notation
\begin{eqnarray*}
 b_k &=& f^{(n-1)}_k\quad ({\rm including}\ b_1 = \omega^{n-1})\\
\bigl[Q(z)\bigr]_k &=& {\rm coefficient \ of} \ z^k\ {\rm in\ the\
polynomial\ } Q(z)
\end{eqnarray*}
Note that \ $a_k\biggl[\bigl(b_1z\bigr)^k\biggr]_k\, =\,
a_k\omega^{(n-1)k}$ but that, by induction,  \ $b_i = f^{(n-1)}_i$ has no
$a_k$'s \  in it if\ $1 \leq i < k$. Thus
we have
\begin{eqnarray*}
f^{(n)}(z) \, &=& \, f\bigr(f^{(n-1)}(z)\bigl) = f\biggl(\sum_{k=1}^{\infty}f^{(n-1)}_kz^k\biggr)
= f\biggl(\sum_{k=1}^{\infty}b_kz^k\biggr)\\
&=& \omega\bigl(b_1z + b_2z^2 + \cdots\bigr) +\ \cdots +\
a_j\bigl(b_1z + b_2z^2 + \cdots\bigr)^j + \cdots \ a_k\bigl(b_1z +
b_2z^2 + \cdots\bigr)^k \ + \ \cdots
\end{eqnarray*}
Considering how $z^k$ can occur in the power\ $\bigl(b_1z + b_2z^2 + \cdots\bigr)^j$, we get
\begin{eqnarray*}
f^{(n)}_k &=& \omega b_k +\ \cdots\ +
 \underbrace{a_j\biggl[\bigl(b_1z +  \cdots + b_{k-j+1}z^{k-j+1}\bigr)^j\biggr]_k}
\ + \cdots\ +
a_k\biggl[\bigl(b_1z\bigr)^k\biggr]_k \\
&\mbox{}&\hskip1.5in\hbox to 1in {{\bf no }\ $a_k$\ {\bf in\ here }}\\
&=& \omega\biggl(a_k\omega^{n-2}%
\bigr(1 + \omega^{k-1} +\cdots +\omega^{(k-1)(n-2)}\bigl)
 + \omega P^{(n-1)}_k(a_1\ldots a_{k-1})\biggr)
+ \cdots \\
&\mbox{}&\qquad\hbox to 2in{} a_j\bigl({\rm polynomial\ with\ no}\ a_k{'s}\bigr) +\ \cdots + a_k\omega^{(n-1)k}\\
  &=& a_k\omega^{n - 1}\bigg(1 + \omega^{k - 1} +
\cdots\  \omega^{(k-1)(n-1)}\bigg) +  P^{(n)}_k(a_1, \ldots, a_{k-1})
\end{eqnarray*}
where $P^{(n)}_k(a_1, \ldots, a_{k-1})$ is the result of summing all the terms
which do not contain $a_k$.
\medskip
By induction on $P^{(n-1)}_k$ and definition of $P^{(n)}_k$, it follows that if \ $a_j=0$ for all $j$ with $1 < j < k$ then \ $P^{(n)}_k(a_1,0, \ldots, 0) = 0$. Moreover, since all computations involve taking integral powers, the coefficients
of $P^{(n)}_k$ are integers. Since the computational process
does not depend on the values of the $a_j$, the polynomial is
independent of $f(z)$.\qquad $\Box$
\bigskip
\section{Proof of  The Conjugation Theorem}
\medskip
{\regsmallcaps Theorem 1A.\ \  [\rm Classification of finite order elements up to conjugation\ ]}

\medskip
 {\sl Suppose that \ $f(z) = \omega z + \sum_{k=2}^{\infty}a_kz^k$ is an
element
of  $G$ of order $n$.\ \  Define $f^\ast$ by
$$f^\ast = \frac{1}{n}\,\sum_{j=1}^n\, \omega^{n-j}f^{(j)}.$$
 Then $f^\ast \in G,\ \ n$ is the multiplicative order of\  $\omega \in F -\{0\}$\
and
$$ f^\ast\circ f\circ\overline{f^\ast}\,(z) = \omega z.$$
Two elements of order $n$ in $G$ are conjugate iff their
lead coefficients are the same primitive $n'$th root of unity.}

\medskip


{\bf Proof:} By Lemma 2.3  the order of $\omega$ is $n$.   From the definition
of composition we have, for all $g, h, k \in
G$:
\begin{eqnarray*}
\bullet \quad \ \ (g + h)\circ k &=& g\circ k \, + \,  h\circ k\\
\bullet \ \quad \ell_\omega \circ(g + h) &=& \ell_\omega \circ g\,  + \, \ell_\omega
\circ h
\end{eqnarray*}
Then we have
\begin{eqnarray*}
\ell_\omega\circ f^\ast & =  \frac{1}{n}\,\sum_{j=1}^n\,
\omega^{n-j+1}f^{(j)},\quad {\rm and}\\
f^\ast \circ f & =  \frac{1}{n}\,\sum_{j=1}^n\,
\omega^{n-j}f^{(j+1)}
\end{eqnarray*}
Since $\omega^n = 1$ and $f^{(n)} = id$, these sums run through the
same terms of $G$. Hence
$$ f^\ast \circ f   =  f_\omega\circ f^\ast \quad
{\rm so\ that}\quad   f^\ast\circ f\circ\overline{f^\ast} = f_\omega. $$
Finally we note from Lemma 2.1 that  any conjugate\
$g\circ f_\omega\circ\overline{g}\,(z) = \omega z\  +\ {\rm higher\ terms}$.
Hence \ $f_\omega$\ is conjugate in $G_1$ to $f_{\omega'}$ if and
only if \ $\omega = \omega'$. So elements of finite order in $G$
are conjugate iff their lead terms are the same. \mbox{$\Box$}
 
\bs 
{\regsmallcaps Theorem 1B.\ \  [\rm Determination of the Conjugating Elements\ ]}

\begin{notation}
$Z_\omega \equaldef\{h \in G\ \vert\  h\circ \ell_\omega
  =\ell_\omega\circ h\} =\ $ the {\bf centralizer} of $\ell_\omega$ in $G$.
\end{notation}

\medskip
H. Furstenberg has pointed out the following fact:

\begin{lemma}
If $\omega$ is a primitive n'th root of unity and $h=h(z) \in G[[z]]$\
then\\  $h \in Z_\omega\  \iff\  h(z) =
\sum_{j=0}^\infty h_{nj+1}z^{nj+1}.$
\end{lemma}
{\bf Proof.}

$ h\circ \ell_\omega(z) =  \sum_{k=1}^\infty h_k(\omega
  z)^k$
\quad  and \quad
$\ell_\omega\circ h(z) = \sum_{k=1}^\infty \omega h_kz^n$.

\medskip
Therefore these are equal \ $ \iff \omega h_k = \omega^k h_k \ \mbox{
 for all}\
 k\ \iff\
h_k = 0$\  for all \ $k\not\equiv 1\,(n).\quad \Box.$

\bigskip
We have the following Corollary to Theorem 1.

\medskip
\begin{corollary}
Suppose that $f(z) = (\omega z + a_2z^2 + a_3z^3  \cdots\ )$ has finite order $n$.
\begin{enumerate}
\item If $g \in G$ then
\[\qquad g\circ f\circ \overline{g} = \ell_\omega\  \iff \ \exists
 h(z) = \sum_{j=0}^\infty h_{nj+1}z^{nj+1}\ \in G\
 \mbox{such that} \ g = h\circ f^\ast.\]
\item For every sequence \ $\{\,g_{nj+1}\}_{j=0}^{\infty}$\
 of
elements of $F$ with $g_1\neq 0$ there exists a unique sequence \ $\{g_k\}_{1 < k \neq nj + 1}$ \  such that
the series $g(z) = \sum_{k=1}^\infty g_kz^k$ satisfies
\ $g\circ f\circ \overline{g} = \ell_\omega$.
\end{enumerate}
\end{corollary}
\medskip
{\bf Proof:}\quad  From Theorem 1,
\[g\circ f\circ \overline{g} = \ell_\omega\  \iff \
g\circ(\overline{f^\ast}\circ \ell_\omega\circ f^\ast)\circ\overline{g} = \ell_\omega \ \iff\ g\circ \overline{f^\ast} \in Z_\omega.\]
Statement (1) now follows from Lemma 3.2.

\bigskip
We use (1) to prove (2):\ A sequence  \ $\{\,g_{nj+1}\}_{j=0}^{\infty}$,\ with $g_1 \neq 0$, extends to a sequence $\{g_k\}_{k=1}$ for which $g(z) = \sum_{k=1}^\infty g_kz^k$\ conjugates $f$ to $\ell_\omega$
iff $g = h\circ f^\ast$ as in (1). We can solve for the missing coefficients of $g$ and the coefficients of $h$ recursively. If\  $f^\ast(z) = b_1z + b_2z^2 + b_3z^3 + \cdots$\
($b_1$ necessarily equals 1) \ then one equates the coefficients of the $z^k$\ in
\[g_1z + g_2z^2 + \ldots g_{n+1}z^{n+1} + \cdots
\ = h_1(z + b_2z^2 + b_3z^3 + \cdots) +
h_{n+1}(z + b_2z^2 + b_3z^3 + \cdots)^{n+1} + \cdots \ ,\]
We leave the details to the reader. \qquad  $\Box$

\bigskip
\section{Proof of the the Construction Theorem for Elements of Order $n$}
{\regsmallcaps Theorem 2:}\quad {[\rm Determination of the coefficients of a formal power series of \\
\hbox{} \hskip .9in finite compositional order]}

\bigskip
 {\sl If\ $n$\ is a positive integer and $\omega$ is a primitive $n'{\rm
th}$\ root of unity in the field $F$ then for every infinite sequence  \
$\{a_k\}_{1 < k \neq nj + 1}$ \ of elements of $F$there exists a unique sequence \ $\{a_{nj+1}\}_{j=1}^{\infty}$
\ such that the formal power series
$$f(z) = \omega z + \sum_{k=2}^{\infty}\ a_k\,z^k$$
 has order $n$ in $G$.}

\bigskip
{\bf Proof of existence:} \quad
 Since $\ell_\omega$ has order
$n$ in $G$, any conjugate \ $\overline{h}\ell_\omega h$\  has order $n$. Given a sequence
\ $\{a_k\}_{1 < k \neq nj + 1}$, \ we shall recursively construct a series $h(z) =
\ \sum_{k=1}^{\infty}\, h_kz^k$\ \ and a sequence $\{a_{nj + 1}\}\ \bigg\vert\ j = 1, 2, \ldots\}$
such that \ $\overline{h}\ell_\omega h(z)$ has the form of the desired $f(z)$.  Indeed we construct $h$ and $f$ simultaneously so that \\
$(h\circ f)(z) = (\ell_\omega\circ h)(z)$. \qquad This is true precisely when
\begin{eqnarray*}(*)\ &\mbox{}& h_1\biggl(\omega z + a_2z^2 + \cdots
 \biggr)  + h_2\biggl(\omega z + a_2z^2 + \cdots \biggr)^2 + \cdots
+  h_k\biggl(\omega z + a_2z^2 + \cdots \biggr)^k + \cdots \\
 = &\mbox{}&\omega h_1z + \omega h_2z^2 + \, \cdots\, \omega h_kz^k +\, \cdots
\end{eqnarray*}
Set $h_1 = 1$. Suppose that \ ${\,h_m}\ (m < k)$\ and \ $a_{nj +1}\ (nj+1 < k)$ have been chosen
so that the coefficients of $z^m \, (m < k)$ on the two sides of  equation (*) are equal.
 Consider the coefficient of \ $z^k$\ on both sides of the above equation:
$$(**)\qquad h_1a_k + \ Q_k\big(h_i, a_\ell \vert \, 2 \leq i < k, \ell < k\big)\  + h_k\omega^k = \omega\,h_k,\quad $$
where \ $h_1 = 1$\ and\ $Q_k$ \ is a polynomial in the given variables.

\medskip
CASE 1:\ \ $k \not\equiv 1(n)$.\quad Note that $a_k$ has been
prespecified and
$\omega^k \neq \omega$. Then we choose according to (**) (and have no choice in so doing):
$$h_k = \frac{a_k + Q_k}{\omega - \omega^k}.$$
CASE 2:\ \ $k \equiv 1(n)$.\quad Note that \ $\omega^k = \omega,\ h_1
= 1$
and to satisfy (**) we choose
\ $$a_k = -Q_k = -Q_k(h_2, \ldots h_{k-1}, \, a_1, \ldots a_{k-1})$$
We may choose \ $h_k$\ arbitrarily. \mbox{$\Box$}

\bigskip
{\bf Remark:}

\begin{itemize}
\item In the above proof, the value of
$a_k$ which we are forced to choose when \ $k \equiv 1(n)$ seems to
depend  on the preceding freely chosen
$h_{nj+1}$.  In fact,we show that Lemma 2.4 implies that  each $a_{nj
  + 1}$
is uniquely determined by the preceding $a_k$'s:
\end{itemize}

{\bf Proof of uniqueness:}
\medskip
If \ an infinite sequence  \
$\{a_k\}_{1 < k \neq nj + 1}$\ is given\ and $f(z) = \omega z +
 \sum_{k=2}^{\infty}\ a_k\,z^k$\quad
 has order $n$ in $G$  then
 \medskip
\begin{itemize}
 \item$f^{(n)}_k = 0$\ for $k > 1$.

 \item $k \equiv 1 ({\rm mod}\ n)\ \implies\ \omega^{k-1} = 1.$
\end{itemize}
 \bigskip
  Hence, for \ $k > 1$\ and\  $k \equiv 1 ({\rm mod}\ n)$,  Lemma 2.4 \ gives \\
$0 = a_k\omega^{n-1}n + P^{(n)}_k(a_1, \ldots, a_{k-1})$.\ Therefore
$$ a_k = -\frac{\omega}{n}P^{(n)}_k(a_1, \ldots,
  a_{k-1}).\qquad\mbox{$\Box$}$$
\bigskip
\section{Further Comments}

 \subsection{On formal power series of order two and of infinite order}\mbox{}

\medskip
From Lemma 2.3 we see that, when $F = {\mathbb R}$, a  real power series of finite order can
only have order one or two and that, for any field of characteristic
0, a series in $G$ of order two is of the form \ $f(z) = -z + a_2z^2 +
a_3z^3 + \cdots$. Note that if $f(z)$ has order two then
$f^\ast(z) =\frac{1}{2}\left(z - f(z)\right)$.

\medskip
L. Shapiro has pointed out that the following is an example of Theorem 1.

\begin{example}[Stanley [14], page 50, problem 41a]\mbox{}

\medskip
Suppose that $f(z) \in G$.

  Then
 $f(-f(-z)) = z \iff
  \exists \ \ g(z) = z + \sum_{n=2}^\infty g_nz^n$ such  that
   $f(z) =\overline{g}(-g(-z))$

\medskip
{\rm In our terms this just says that  $f\circ \ell_{-1}$ has
order two iff it is conjugate to $\ell_{-1}$}
\end{example}

\bigskip

\subsection{Muckenhoupt's results on formal series of infinite order}\mbox{}

\medskip
From Lemma 2.3 we see that \ $f(z) = (\omega z + a_2z^2 + \cdots$\ )
has infinite order in $G$ if $\omega$ has infinite multiplicative
order in $F$. And from Theorem 2 we see that if $\omega$ has finite
order, it is still possible that $f(z)$ has infinite order -- one
can choose any one of the $a_{nj+1}$ in such a way that $f(z)$
won't have finite order.  In these two situations, Muckenhoupt [8]
has explained what happens (while implying that the following results were classically known):
 \medskip
\begin{unnumberedtheorem}[Muckenhoupt]
Suppose that $f(z) = (\omega z + a_2z^2 + a_3z^3 + ...) \in G$
\begin{enumerate}
\item If $\omega$ has infinite order  then $f$ is conjugate to
$\ell_\omega$.
\item If $\omega$ has finite order $n$ then $f$ is conjugate to
\ $g \in G$, where $g$ is of the form

$g(z) = \omega z + \sum_{j=1}^\infty g_{nj+1}z^{nj+1}$.
\end{enumerate}
\end{unnumberedtheorem}
\medskip
{\bf Remark:} Our first proof that formal power series of finite
order under composition are conjugate
to their linear parts used an inductive argument
from Lemmas 2.3 and 2.4 above and Part (2) of Muckenhoupt's theorem, without explaining the role of $f^\ast$.
   
\section{Useful References and Comments on the History of Formal Power Series}
\ms
The subject of Formal Power Series is a deep classical subject of importance in Analysis, Combinatorics and Algebra.  I hope the comments below will be helpful and will add some perspective, but there is no claim to completeness and I apologize to those contributors to the subject whom I omit. 
\ms
\subsection{Introductory surveys}
\ms
Excellent introductions to formal power series and surveys of the field are given by
\be
\item Ivan Niven [9], (1969) -- this article won the MAA's Lester R. Ford Award for expository excellence.
\item Henri Cartan [4], Section I.1 -- originally published 1961.
\item Stephen Scheinberg [12], (1970)
\item The volumes [15], [16] of Richard P. Stanley (1997, 1999).
\item Thomas S. Brewer [2], (2014).
\item Anthony G. Farrell and Ian Short [10], Ch. 10 (2015).
\ee
\subsection{ E.  {Schr\"{o}der}'s  foundational 1871 paper } \mbox{}

\ms
Ernst \schroeder's  foundational 1871 paper [13], {\it ``Ueber Iterate Functionen"} (``On Iterated Functions")  set the stage for much to come. In his opening he states (translated)

\ms
\bc
\parbox{5in}{\em `` I present here investigations into a field in which I have seen very\\ \mbox{ \ \ } little previous work."}
\ec

\ms
This rich paper formally introduces the concept of the $r^{\rm th}$ iterate of a function $F(z)$ and, in particular, contains the following two items which relate directly  to the work above; namely, to Lemma 2.4 on the determination of $f_k^{(n)}$ and to the conjugation problem in Theorem 1.

\ms
\be
\item \S6, page 310:  {\bf Coefficients of $F^{(r)}(z)$ when $F(z) = z + a_2z^2 + \ldots$ is a MacLaurin series.}

\ms
\noindent Between pages 310 and 315, this gives a remarkably deep and detailed calculation of the $k^{\rm th}$ coefficient of $F_k^{(r)}$ (which \schroeder \ denotes $F_k^r(z)$). Since he does not consider convergence questions until after he develops these formulas, the development of his formulas really takes place within the realm of formal power series over $\mathbb{C}$.  I conjecture that there is much to be gained by those who can wield his formulas in particular cases.

\ms
\item \S3(D), page 303:\ \ Schr\"{o}der assumes given a function $F(\zeta)$ of the variable $\zeta$.  He points out how useful it would be in evaluating the $r'th$ iterate $F(F( \cdots(\zeta)\cdots ))$, if there existed a constant $m$ and a change-of-coordinates function $\psi$ such that

\ms
\bc
\mbox{\hspace{1in}}\framebox{\parbox{1in}{\centerline{$\psi(m\zeta) = F{\psi(\zeta)}$}}}\quad {(\bf ``{Schr\"{o}der's} Equation'')}
\ec
\ms
In our notation, if $\psi$ is an invertible function, we would then have \break $F= \psi\circ\ell_m\circ\psi^{-1}$ and it becomes possible to easily compute
\[F^{(r)}(\zeta)=\psi(m^r\cdot\psi^{-1}(\zeta)).\]

\noindent
The solution of {Schr\"{o}der's} 
equation has been sought in many different settings and formulation of the Conjugacy Problem is often phrased in terms of ``solving \schroeder's equation".
\ee

\subsection{Conjugacy of Formal Power Series}
\be
\item {\bf Great Historical Papers Can Be Used Incorrectly.}

\ms
\noindent
The algebraic subject of Formal Power Series over $\mathbb{C}$ is so deeply intertwined with the subject of MacLaurin series in Analysis, that one can easily misapply a historical result from one  of these subjects incorrectly in  the other subject. 

\ms
For example, the following mistakes are possible:

\ms
\noindent
\be
\item Suppose given a MacLaurin series  $F(z) = a_1z + a_2z^2 + \ldots$ with positive radius of convergence. Seeking an analytic function which linearizes a curvilinear angle,  one might formally find a conjugator $\psi(\zeta)$ satisfying {Schr\"{o}der's} equation which is not analytic (does not converge in any neighborhood of the origin) and mistakenly think that they've found an analytic solution. In 1917, Pfeiffer [11] pointed out this problem by giving an analytic $F(z)$, where $a_1$ has infinite order, such that there exists no analytic solution to {Schr\"{o}der's} equation.

\ms
It is noteworthy that Pfeiffer does prove (page 186 of [11]) that, when a formal power series $F(z)$ has $a_1$ of infinite order, a unique {\bf\em formal} solution to {Schr\"{o}der's} equation always exists. 
\ms
\item A person wishing to prove Theorem 1A, given a formal power series of finite order, $F(z) = a_1z + a_2z^2 + \ldots$,  might wish to use a paper in analysis -- such as  Siegel's famous paper [14]-- to  show that this can be conjugated to a linear function.  {\bf However}, Siegel's beautiful argument on the opening page of [14] assumes that the function being conjugated is analytic in a neighborhood of the origin and is stable (a valid assumption for analytic functions of finite compositional order). He then finds the conjugator with a mixture of topology (the universal covering space) and analysis (Schwartz' lemma). However, this analyticity assumption is invalid in proving Theorem 1A, since we have proven in the Corollary to Theorem 2 that not every formal power series of finite order is analytic.
\ee

\ms
\noindent{\bf But note:} Despite these warnings about the possible misuse of the historical documents, 
\bc
{\bf \em the beauty and insights of the early papers\\ is not to be bypassed. }
\ec 
And often it happens that, inside an argument in one field there will be a gem which can be of great use in the other.  For example,  on page 42, line 3, of [3], in the midst of an analysis argument which assumes that ``{\it the sequence of iterates $\{f^n\}$ is uniformly bounded in some neighborhood of the origin}" the formula for the sum which we call $f^\ast$ in Theorem 1A appears.  If $f$ is an analytic function of finite order, the conjugating function $\varphi$ which that proof gives is in fact $\varphi = f^\ast$.  Seeing this one could conjecture that $f^\ast$ is also the conjugating element in the group of formal power series, and come up with an algebraic proof. 

\ms
\noindent One cannot help but conjecture that the great mathematicians of the past, having the proof of the conjugacy-to-linear-term in the analytic setting, knew from this a proof in the formal power series setting.  Indeed, the 1948 proof of Bochner and Martin [1] mentioned above notedly avoids compactness considerations in the complex plane by putting ``the ordinary weak topology based on convergence for each coefficient separately" on the set of formal power series directly, and proving a theorem (Theorem 8) about {\bf \em bounded} groups of formal power series.  Clearly any formal power series of finite order generates a finite, hence bounded, group of formal power series.
\ms
\item {\bf Scheinberg's Table of Normal Forms of Conjugation}

\ms
\noindent Scheinberg [12] gives a very rich discussion  of formal power series as of 1970, with interesting historical detail.  He includes at the end of \S3 with Table 1, a table of canonical forms of formal power series under conjugation, including, for example, Muckenhoupt's Theorem above.

\ms
\noindent
\item {\bf Conjugating Formal Power Series of Finite Order}

\ms
\noindent Summarizing much of what has been said above, the fact (Theorem 1 above) that a formal power series $f(z) = \omega z + \sum_{k=2}^{\infty}a_kz^k$ of compositional order $n$ is conjugate in the group $G$ of formal power series to $\ell_\omega(z) = \omega z$ is very well known.  It has appeared at least in [1], [2], [5], [10] and [12]. 

\ms
\noindent The fact that $f^\ast$ of Theorem 1 may be used as the conjugator appears in[ 1], [5] and [10], and has been known at least since Bochner and Martin's 1948 book [1].

\ms
\noindent {\bf Historical challenge:} 

\ms
 Check out all of the references given in Peiffer's 1917 paper [11] and see whether the knowledge of this result goes back to the early twentieth or even the nineteenth century.
\ee

\ms
\noindent
\subsection{The Construction of Formal Power Series of Finite Order}

\mbox{}

\ms
\noindent Independent proofs of Theorem 2 above were given by the author and Thomas S. Brewer [2]. We do not know (May, 2018) of other proofs, but in this venerable field, we would hesitate to  make brash claims.

\bigskip
\hrulefill

\bigskip

\centerline{\large REFERENCES}
\bigskip
\begin{enumerate}
\item {\bf Salomon Bochner and William Ted Martin}, {\it Several Complex
Variables}, Princeton University Press (1948).
\item {\bf Thomas S. Brewer}, {\it Algebraic Properties of Formal Power Series Composition}, University of Kentucky Dissertation (2014).
\item{\bf Lennart Carleson and Theodore W. Gameli}, {\it Complex Dynamics}, Springer-Verlag New York, Inc., Universitext: Tracts in Mathematics (1993).
\item {\bf Henri Cartan}, {\it Elementary Theory of Analytic Functions of One Or Several Variables},
Dover Books (1995).
\item {\bf Gi-Sang Cheon and Hana Kim}, {\it The Elements of finite order in the Riordan group over the complex field}, Linear Algebra and Its Applications 439(2013), 4032 - 4046.
\item {\bf Marshall M. Cohen}, {\it Elements of Finite Order in the Group of Formal Power Series Under Composition}, MAA Section Meeting April 8, 2006.
\item {\bf S. A. Jennings}, {\it Substitution Groups of Formal Power Series}, Canadian J. Math 6 (1954), 325-340.
\item {\bf Benjamin Muckenhoupt}, {\it Automorphisms of Formal Power Series Under \break
 Substitution}, Transactions
Amer. Math. Soc. 99 (1961), 373--383.
\item  {\bf Ivan Niven},  {\it Formal Power Series}, Amer. Math. Monthly 76 (1969), 871 - 889.
\item {\bf Anthony G. O'Farrell and Ian Short}, {\it Reversibility in Dynamics and Group Theory}, London Mathematical Society Lecture Note Series: 416, Cambridge University Press (2015)
\item {\bf G. A. Pfeiffer}, {\it On the Conformal Mapping of Curvilinear Angles.  The Functional Equation 
$\phi[f(x)] = a_1\phi(x)$}, Trans. Amer. Math. Soc. 18 (1917), p. 185 -198.
\item {\bf Stephen Scheinberg}, {\it Power Series in One Variable}, Journal of Mathematical Analysis and Applications 31,  p. 321-333 (1970).
\item {\bf Ernst Schr\"{o}der}, {\it Ueber Iterirte Functionen}, Math. Annal. 3 (1871), p. 296 - 322.
\item {\bf Carl Ludwig Siegel}, {\it Iteration of Analytic Functions}, Annals of Mathematics, Vol. 43, No. 4, (October, 1942), p. 607 - 612.
\item {\bf Richard P. Stanley}, {\it Enumerative Combinatorics, Volume 1, Section 1.1}, Cambridge Studies in Advanced Mathematics 49, Cambridge University Press (1997).
\item  {\bf Richard P. Stanley}, {\it Enumerative Combinatorics, Volume 2, Section 5.4}, Cambridge Studies in Advanced Mathematics 62, Cambridge University Press (1999).
\end{enumerate}

\hrulefill
\medskip

\end{document}